\documentclass[a4paper,12pt]{article}

\newcommand{\E}  {\mathbb{E}}
\renewcommand{\P}  {\mathbb{P}}

\usepackage{rotating}
\usepackage{amsmath, amssymb, amsfonts}
\newcommand\unnumberedfootnote[1]{ %
        \let\temp=\thefootnote %
        \renewcommand{\thefootnote}{}%
        \footnote{#1}%
        \let\thefootnote=\temp%
        \addtocounter{footnote}{-1}}

\newcommand{\N}{\mathbb{N}}
\newcommand{\Z}{\mathbb{Z}}
\newcommand{\R}{\mathbb{R}}
\renewcommand{\P}{\mathbb{P}}
\renewcommand{\L}{\mathcal{L}}

\def \bpf {\noindent{\sc Proof: }}
\def \epf {\hbox{}\nobreak\hfill
\vrule width 2mm height 2mm depth 0mm
\par \goodbreak \smallskip}

\usepackage{graphicx}
\usepackage{graphics}
\usepackage{pstricks}
\usepackage{subfigure}
\usepackage{color}
\definecolor{DarkGreen}{rgb}{0,0.6,0}
\definecolor{VertFonce}{rgb}{0,0.6,0}
\definecolor{MidGreen}{rgb}{0.6,1,0.6}
\definecolor{LightGreen}{rgb}{0.88,1,0.88}
\definecolor{LightGray}{rgb}{0.94,0.94,0.94}
\definecolor{VeryLightBlue}{rgb}{0.9,0.9,1}
\definecolor{LightBlue}{rgb}{0.8,0.8,1}
\definecolor{DarkBlue}{rgb}{0,0,0.6}
\definecolor{VeryLightYellow}{rgb}{1,1,0.9}
\definecolor{LightYellow}{rgb}{1,1,0.6}
\definecolor{MidYellow}{rgb}{1,1,0.5}
\definecolor{VeryLightRed}{rgb}{1,0.9,0.9}
\definecolor{LightRed}{rgb}{1,0.8,0.8}
\definecolor{turquoise}{rgb}{0.00,0.53,0.68}{}
\definecolor{mauve}{rgb}{0.50,0.00,0.50}

\newtheorem{theorem}{Theorem}[section]
\newtheorem{lemma}[theorem]{Lemma}
\newtheorem{proposition}[theorem]{Proposition}
\newtheorem{corollary}[theorem]{Corollary}

\newtheorem{remark}[theorem]{Remark}
\numberwithin{equation}{section}

\numberwithin{equation}{section}
\begin{document}
\title{\LARGE{Branching processes with competition and generalized Ray Knight Theorem}
\author{Mamadou Ba \and Etienne Pardoux  }}
\maketitle \noindent{\footnotesize LATP-UMR 6632, CMI, Universit\'e de Provence, 39 rue F. Joliot Curie,\\
 Marseille cedex 13, FRANCE.  \\
{\tt email : mba@cmi.univ-mrs.fr ; pardoux@cmi.univ-mrs.fr}
}
~~\\~~\\~~\\
\begin{abstract}
{\small 
We consider a discrete model of population with interaction  where the  birth and death  rates  are  non linear functions of the population size. 
After proceeding to renormalization of the model parameters, we obtain in the limit of large population that the population size evolves as a diffusion solution of the SDE
\begin{align*} 
 Z^x_t =x+\int_0^t f(Z^x_s)ds+2\int_0^t\int_0^{Z^x_s}W(ds,du),
\end{align*}
where $W(ds,du)$ is a time space white noise on $([0,\infty))^2$.\\
We give a Ray-Knight  representation of this diffusion in terms of the local times of a reflected Brownian motion $H$ with a drift that depends upon the local time accumulated by $H$ at its current level, through the function $f'/2$.
}
 \end{abstract}
\unnumberedfootnote{{\sc Keywords} : Galton--Watson processes with interaction, generalized Feller diffusion}.
\unnumberedfootnote{\emph{AMS 2000 subject classification.}
  (Primary) {\tt 60J80, 60F17} (Secondary) {\tt 92D25}.}
  \newpage
\section*{Introduction}
Consider a population evolving in  continuous time with $m$ ancestors at time $t=0$, in which each individual, independently of the others,  gives birth to children at a constant rate $\mu$, and dies  after an exponential time with parameter $\lambda$.  For each individual we superimpose additional birth and death rates due to interactions  with others  at a certain rate which depends upon the size of the total population. For instance, we might decide that each individual  dies because of competition at a rate equal to $\gamma$ times the number of presently alive individuals in the population, which amounts to add a global death rate equal to $\gamma (X^m_t)^2$, if $X^m_t$ denotes the total number of alive individuals at time $t$. 
% It is rather clear that the process which describes the evolution of the total population, which is not a branching process (due to the interactions between branches, created by the competition term), goes 
%extinct in finite time a. s.

If we consider this population with $m=[Nx]$ ancestors at time $t=0$, weight each individual with the factor $1/N$, and choose
$\mu_N=2N+\theta$, $\lambda_N=2N$ and $\gamma_N=\gamma/N$, then it is shown in Le, Pardoux and Wakolbinger \cite{4LPW} in the above particular case of a quadratic competition term
that the ``total population mass process'' converges weakly to the solution of the Feller SDE with logistic drift
\begin{align}\label{4eqZ}
 dZ^x_t=\left[ \theta Z^x_t- \gamma (Z^x_t)^2\right]  dt+2\sqrt{Z^x_t}dW_t,~ Z^x_0=x.
\end{align}
The diffusion $Z^x$ is called Feller diffusion with logistic growth and  models the evolution of the size of a large population with competition. In this model
$\theta$ represent the supercritical branching parameter while  $\gamma$ is the rate at which each  individual is killed by any one of his contemporaneans.  This model has been studied in Lambert \cite{4LA}, who shows in particular that its extinction time is finite almost surely. 

We generalize the logistic model by replacing the quadratic function $\theta z-\gamma z^2$  by a more general nonlinear function $f$ of the population size. We then obtain in the continuous setting  a diffusion which is the solution of the SDE
\begin{align}\label{4fA}
 Z^x_t &=x+\int_0^t f(Z^x_s)ds+2\int_0^t\int_0^{Z^x_s}W(ds,du),
\end{align}
where the function $f$ satisfies the following hypothesis.\\
 \textbf{Hypothesis A:} $f\in C(\R_+;\R)$, $f(0)=0$ and $\exists \beta\ge 0$ such that 
\begin{align*} 
f(x+y)-f(x)\le \beta y\quad \forall x,y\ge 0.
\end{align*}
The equation \eqref{4fA} has a unique strong solution (see \cite{4DL}).
Note that the hypothesis \textbf{A} implies that 
\begin{align*}
\forall x\ge 0, f(x)\le \beta x.
\end{align*}
An equivalent way to write \eqref{4fA} is the following.
\begin{align}\label{4fA1}
 Z^x_t =x+\int_0^t f(Z^x_s)ds+2\int_0^t \sqrt{Z^x_s}dW^x_s,
\end{align}
where $W^x$ is a standard Brownian motion.
However, the joint evolution of the various population sizes $\left\{ Z^x_t,t\ge 0\right\}$ corresponding to different initial population sizes $x$ would necessitate a complicated description of the joint law of the $\left\{ W^x_., x\ge 0\right\}$. Whereas the formulation \eqref{4fA} due to Dawson, Li \cite{4DL} with one unique space-time white noise $W$, describes exactly the joint evolution of $\left\lbrace Z^x_t,t\ge 0,x\ge 0\right\rbrace$ which we have in mind. 
%The difference between \eqref{4fA} and \eqref{4fA1} is  that the Brownian motion $W^x$ in \eqref{4fA1} depends on  the %initial condition $x$ instead the time space white noise $W(ds,du)$ does not.
%Existence and uniqueness of a strong solution of \eqref{4fA1} follows from Proposition 3.14 in \cite{4PR}. 
We call this diffusion  the generalized Feller diffusion. In order to derive  this continuous model, we first define a discrete model.  For defining jointly  the discrete model for all initial population sizes, we need as in \cite{4PW} to impose a non symmetric competition rule between the individuals, which we will describe in section 1 below. We do a suitable renormalization of the parameters of the discrete model in order to obtain  in section 2 a large population limit of our model which is a generalized Feller diffusion. Section 3 is devoted to give  a Ray Knight representation for such a generalized Feller diffusion. The proof of this representation uses tools from stochastic analysis, in particular the ``excursion filtration", following an analogous proof of another generalized Ray Knight theorem in \cite{4NRW}. 

\section{Discrete model with a general interaction}\label{4discrete}
In this section we set up a discrete  mass continuous time approximation of the generalized Feller diffusion. We consider a discrete model of population with interaction  in which each individual, independently of the others, gives birth naturally at rate $\lambda$, dies naturally at rate $\mu$. Moreover, we suppose that each individual gives birth and dies because of  interaction  with others at  rates which depend upon  the current population size. Moreover, we exclude multiple births at any given time and we define the interaction  rule through  a function $f$ which satisfies hypothesis \textbf{A}.

In order to define  our model jointly for all initial sizes, we need to introduce a non symmetric description of the effect of the interaction as in \cite{4BP} and \cite{4LPW}, but here we allow the interaction to be favorable to some individuals.
\subsection{The model} 
We consider a continuous time $\Z_+$--valued population process
$\{X^m_t,\ t\ge0\}$, which starts at time zero from  $m$ ancestors who are arranged  from left to right, and evolves in continuous time. The left/right order  is passed on to their offsprings: the daughters are placed on the right of their mothers  and  if at a time $t$ the individual $i$ is located at the left of individual $j$, then all the offsprings of $i$ after time $t$ will be placed on the left of all the offsprings of $j$. Since we have excluded multiple births at any given time, this means that the forest of genealogical trees of the population is a planar forest of trees, where the ancestor of the population $X^1_t$ is placed on the far left, the ancestor of $X^2_t-X^1_t$ immediately on his right, etc... Moreover, we draw the genealogical trees in such a way that distinct branches never cross. This defines in a non--ambiguous way an order from left to right within the population alive at each time $t$. See Figure 1.
\begin{figure}[htbp]
 \centering
  \includegraphics[width=0.85\textwidth]{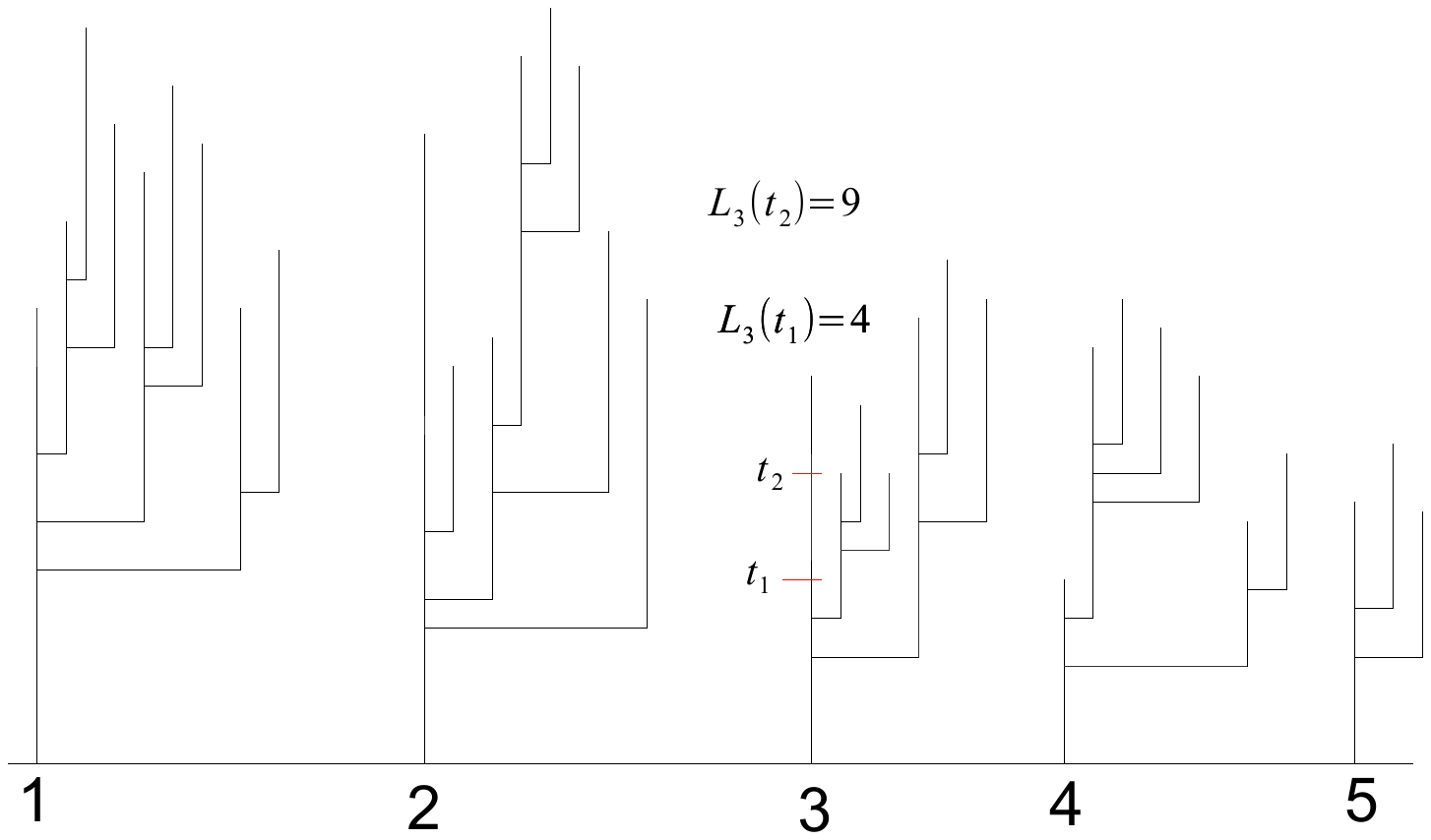} 
 \caption {Planar forest with five ancestors}
 \label{4fig1}
 \end{figure}
 We decree that  each individual feels the interaction with the others placed on his left but  not with those on his right. Precisely, at any time $t$, the individual $i$ has an interaction death rate  equal to $\left( f(\L_i(t)+1)-f(\L_i(t))\right)^-$ or an interaction birth rate equal to $\left( f(\L_i(t)+1)-f(\L_i(t))\right)^+$, where $\L_i(t)$ denotes the number of individuals alive at time $t$ who are located on the left of $i$ in the above planar picture.
This means that the  individual $i$ is under attack by the others located at his left  if  $f(\L_i(t)+1)-f(\L_i(t)) <0$ while the interaction improve his fertility if $f(\L_i(t)+1)-f(\L_i(t))>0$.  
Of course, conditionally upon $\L_i(\cdot)$, the occurence of a ``competition death event''  or an ``interaction birth event"  for individual $i$ is independent  of the other birth/death events and of what happens to the other individuals. In order to simplify our formulas, we suppose moreover that  the first individual in the left/right order has a birth rate equal to $\lambda+f^+(1)$ and a death rate equal to  $\mu+f^-(1)$.

The resulting total interaction death and birth  rates endured by the population $X^m_t$ at time $t$ is then
$$\sum_{k=1}^{X^m_t}[(f(k)-f(k-1))^+ -(f(k)-f(k-1)^-]=\sum_{k=1}^{X^m_t}(f(k)-f(k-1))= f(X^m_t).$$

As a result, $\{X^m_t,\ t\ge0\}$ is a continuous time $\Z_+$--valued Markov process, which evolves as follows. $X^m_0=m$. If $X^m_t=0$, then $X^m_s=0$ for all $s\ge t$. While at state $k\ge1$, the process 
$$X^m_t\text{ jumps to }\begin{cases}k+1,&\text{ at rate $\lambda k+\sum_{\ell=1}^{k}(f(\ell)-f(\ell-1))^+$};\\
                               k-1,&\text{ at rate $\mu k+\sum_{k=1}^{k}(f(\ell)-f(\ell-1))^-$}.
                               \end{cases}  $$
                               
  \subsection{Coupling over ancestral population size}
The above description specifies the joint evolution of all $\{X^m_t,\ t\ge0\}_{m\ge 1}$, or in other words of the two--parameter process $\{X^m_t,\ t\ge0, m\ge1\}$.  
 In the case of a linear function $f$, for each fixed $t>0$, $\{X^m_t,\ m\ge1\}$ is an independent increments process.
 In the case of a nonlinear function $f$, we believe that for $t$ fixed $\{X^m_t,\ m\ge1\}$ is not a Markov chain. 
 That is to say, the conditional law of $X^{n+1}_t$ given $X^n_t$ differs from its conditional law given 
 $(X^1_t,X^2_t,\ldots,X^n_t)$. The intuitive reason for that is that the additional information carried by
 $(X^1_t,X^2_t,\ldots,X^{n-1}_t)$ gives us a clue as to the fertility or the level of competition that  the progeny
 of the $n+1$st ancestor had to beneficit or to suffer from, between time 0 and time $t$.
 
 However, 
 $\{X^m_\cdot,\ m\ge1\}$ is a Markov chain with values in the space $D([0,\infty);\Z_+)$
 of c\`adl\`ag functions from $[0,\infty)$ into $\Z_+$, which starts from 0 at $m=0$. Consequently, in order to describe the law of the whole process, that is of the two--parameter process $\{X^m_t,\ t\ge0, m\ge1\}$, it suffices to describe the conditional law of $X^n_\cdot$,
 given $\{X^{n-1}_\cdot\}$. We now describe that conditional law for arbitrary $1\le m<n$. Let $V^{m,n}_t:=X^n_t-X^m_t$,
 $t\ge0$. Conditionally upon $\{X^\ell_\cdot,\ \ell\le m\}$, and given that $X^m_t=x(t)$, $t\ge0$, $\{V^{m,n}_t,\ t\ge0\}$ is a $\Z_+$--valued
 time inhomogeneous
 Markov process starting from $V^{m,n}_0=n-m$, whose time--dependent infinitesimal generator $\{Q_{k,\ell}(t),\ k,\ell\in\Z_+\}$ is such that its off--diagonal terms are given by
 \begin{align*}
 Q_{0,\ell}(t)&=0,\quad \forall \ell\ge1,\quad \text{ and for any } k\ge1,\\
 Q_{k,k+1}(t)&=\mu k+\sum_{\ell=1}^{k}(f(x(t)+\ell)-f(x(t)+\ell-1))^+, \\
 Q_{k,k-1}(t)&=\lambda k+\sum_{\ell=1}^{k}(f(x(t)+\ell)-f(x(t)+\ell-1))^-,\\
 Q_{k,\ell}(t)&=0,\quad\forall\ell\not\in\{k-1,k,k+1\}.
 \end{align*}

The reader can easily convince himself that this description of the conditional law of 
$\{X^n_t-X^m_t,\ t\ge0\}$, given $X^m_\cdot$ is prescribed by what we have said above, 
and that  $\{X^m_\cdot,\ m\ge 1\}$ is indeed a Markov chain.
\begin{remark}
Note that  if the function $f$ is  increasing on [0, $a$], $a>0$ and decreasing  on $[a,\infty)$, the interaction improves the rate of fertility  in a population whose size is smaller than $a$ but for large size  the interaction amounts to competition  within the population. This  is reasonable because when the population is large,  the limitation of resources implies competition within the population.   For a positive interaction (for moderate population sizes) one  can realize that an increase in the population size allows a more efficient organization of the society, with specalisation among its members, thes resulting in better food production, health care, etc...  We are mainly  interested in the model with  interaction defined with functions $f$ such that $\lim_{x\to\infty}f(x)=-\infty.$ 
Note also that we could have generalized our model to the case $f(0)\ge 0$. $f(0)>0$ would mean an immigration flux. The reader can easily check that results in section 2 would still be valid in this case. However in Proposition \ref{subcritic} and in section 3.3 below, assumption $f(0)=0$ is crucial, since we need the population to get extinct in finite time $a.s.$.
\end{remark}
\subsection{The associated exploration process in the discrete model} 
The just described reproduction dynamics gives rise to a forest $\mathcal{F}^m$ of $m$ trees of descent, drawn into the plane as sketched in Figure 2.   Note also that, with the above described construction, the $\left( \mathcal{F}^m, m\ge 1\right)$, are coupled: the forest $\mathcal{F}^{m+1}$ has the same law as the forest $\mathcal{F}^m$ to which we add a new tree generated by an ancestor placed at the $(m+1)$st position. If the function $f$ tends to $-\infty$ and $m$ is large enough, the trees further to the right of the forest $\mathcal{F}^m$ have a tendency to stay smaller because of the competition : they are ``under attack" from the trees to their left.
From $\mathcal{F}^m$ we read off a continuous and piecewise linear $\R_+ $-valued path $H^{m} = \left( H^{m}_s \right)$ (called the exploration process of $\mathcal{F}^m$) which is described as follows.

Starting from the initial time $s=0$ the process $H^m$ rises  at speed $p$ until it hits the top of the first ancestor branch (this is the leaf marked with $D$ in Figure 2). There it turns and goes downwards, now at speed $-p$, until arriving at the next branch point (which is $B$ in Figure 2). From there it goes upwards into the (yet unexplored) next branch, and proceeds in a similar fashion until being back at height 0, which means that the exploration of the leftmost tree is completed. Then explore the next tree, and so on. See Figure 2.

\begin{figure}[htbp]
\begin{center} 
\includegraphics[width=0.85\textwidth]{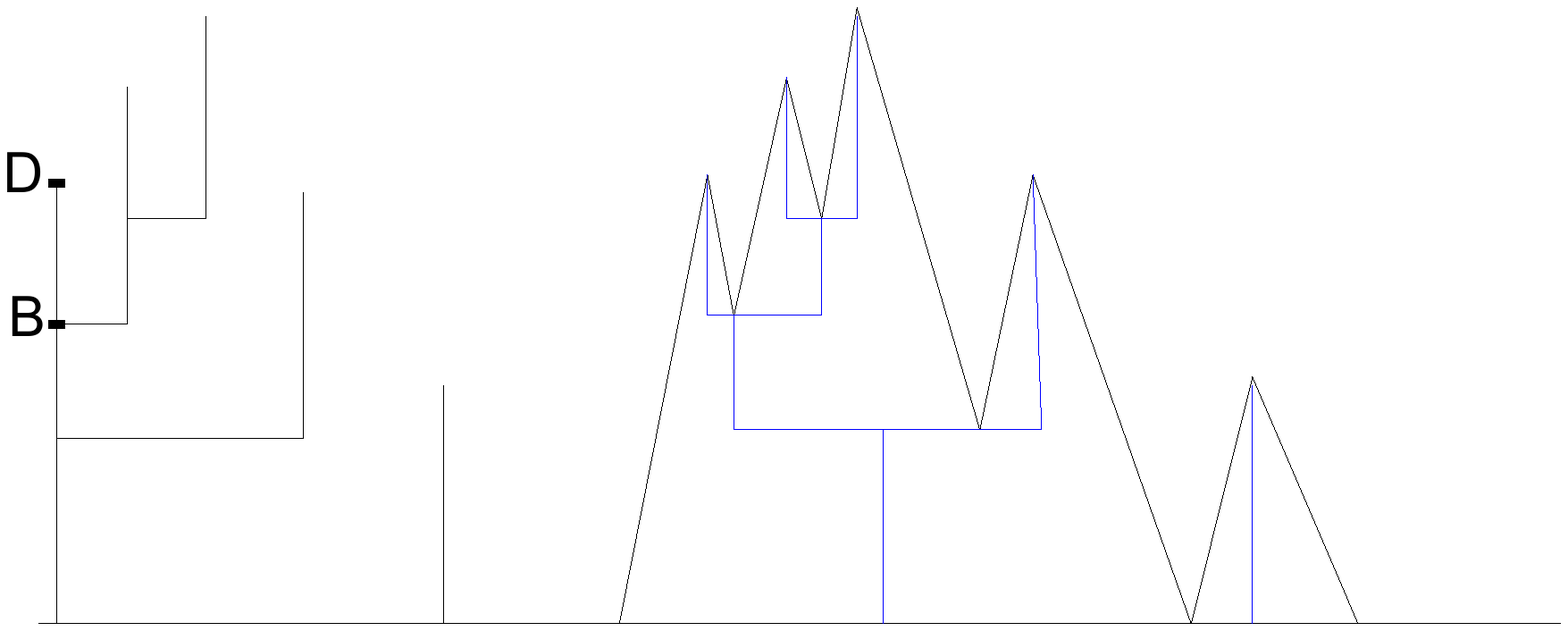}
\caption{A forest with two trees and its exploration process.}
\end{center}
\label{4fig2}
\end{figure}
\begin{figure}[htbp]
\begin{center} 
\includegraphics[width=0.85\textwidth]{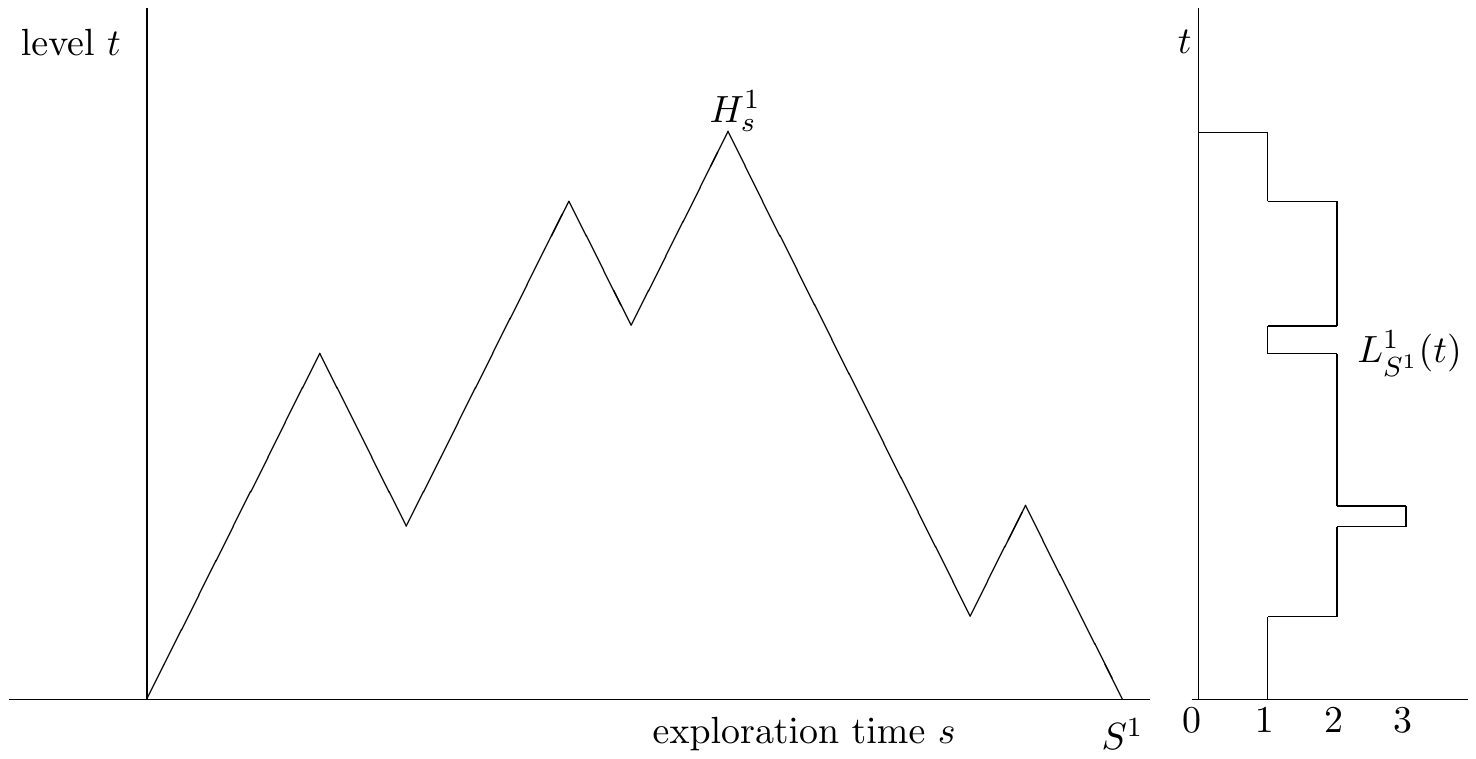}
\caption{Discrete Ray Knight representation.}
\end{center}
\label{4fig4}
\end{figure}

We  define the local time $L^m_s(t)$ accumulated by  the process $H^m$ at level $t$ up to time $s$ by:
\[L^m_s(t)=\lim_{\epsilon\rightarrow 0}\frac{1}{\epsilon}\int_{0}^s 1_{t \leq H^m_r<t+\epsilon} dr.\]
The process $H^m$ is  piecewise linear, continuous with derivative  $\pm p$ : at any time $s\ge 0$, the rate of appearance of minima (giving rise to new branches) is equal $$p\mu+\left[f(\lfloor \frac{p}{2}L^m_s(H^m_s)\rfloor+1)-f(\lfloor \frac{p}{2}L^m_s(H^m_s)\rfloor)\right]^+ ,$$ and the rate of appearance of maxima (describing deaths of branches)  is equal to $$p\lambda +\left[f(\lfloor \frac{p}{2}L^m_s(H^m_s)\rfloor+1)-f(\lfloor \frac{p}{2}L^m_s(H^m_s)\rfloor )\right]^-.$$ Let  $S^m$ be the time needed in order to explore the forest $\mathcal{F}^m$. We have 
$$S^m=\inf\lbrace s>0; \frac{p}{2}L^m_s(0)\ge m \rbrace .$$
Under the assumption that $S^m<\infty$ $a.s.$ for all $m\ge 1$, we have the following discrete Ray Knight representation (see Figure \ref{4fig4}). 
\begin{align*}
\left( X^m_t, t\ge 0, m\ge 1\right)\equiv\left( \frac{p}{2}L^m_{S^m}(t),t\ge 0, m\ge 1\right).
\end{align*}
\subsection{Renormalized discrete model}\label{4secRenorm}
Now we proceed to a renormalization of this model. For $x\in \R_+$ and $N\in\N$, we choose $m=\lfloor Nx\rfloor$, $\mu=2N$, $\lambda=2N$,  we multiply $f$ by N and divide by $N$ the argument of the function $f$. We affect to each individual in the population a mass equal to $1/N$. Then the  total mass process $Z^{N,x}$, which starts from $\frac{\lfloor Nx\rfloor}{N}$ at time $t=0$, is a Markov process whose evolution can be described as follows.\\
$Z^{N,x}$ jumps  from $\frac{k}{N}$ to
$\begin{cases} \frac{k+1}{N} \text{ at rate }  2Nk+N\sum_{i=1}^k\left(f(\frac{i}{N})-f(\frac{i-1}{N})\right)^+\\ 
\frac{k+1}{N} \text{at rate } 2Nk+N\sum_{i=1}^k\left((f(\frac{i}{N})-f(\frac{i-1}{N})\right)^-,
\end{cases}$ \\
Clearly there exist two mutually independent standard Poisson processes $P_1$ and
$P_2$ such that
\begin{align*}
Z^{N,x}_t&=\frac{\lfloor Nx\rfloor}{N}+\frac{1}{N}P_1\left(\int_0^t \left(2N^2 Z^{N,x}_r+N\sum_{i=1}^{NZ^{N,x}_r}(f(\frac{i}{N})-f(\frac{i-1}{N}))^+\right)dr\right)\\
&-\frac{1}{N}P_2\left(\int_0^t\left(2N^2 Z^{N,x}_r+N\sum_{i=1}^{NZ^{N,x}_r}(f(\frac{i}{N})-f(\frac{i-1}{N}))^-\right)dr\right).
\end{align*}
Consequently there exists a local  martingale $M^{N,x}$ such that
\begin{equation}\label{4fAN}
Z^{N,x}_t=\frac{\lfloor Nx\rfloor}{N}+\int_0^t f(Z^{N,x}_r) dr+M^{N,x}_t.
\end{equation}
Since $M^{N,x}$ is a purely discontinous local martingale, its quadratic variation $[M^{N,x}]$ is given by the sum of the squares of its jumps, i.e. 
\begin{align}\label{4crodroit}[M^{N,x}]_t&=\frac{1}{N^2}\left[P_1\left(\int_0^t \left(2N^2 Z^{N,x}_r+N\sum_{i=1}^{NZ^{N,x}_r}(f(\frac{i}{N})-f(\frac{i-1}{N}))^+\right)dr\right)\right. \nonumber\\
&\left.+ P_2\left(\int_0^t\left(2N^2 Z^{N,x}_r+N\sum_{i=1}^{NZ^{N,x}_r}(f(\frac{i}{N})-f(\frac{i-1}{N}))^-\right)dr\right)\right].   \end{align}

We deduce from \eqref{4crodroit} that the conditional quadratic variation $\langle M^{N,x}\rangle$ of $M^{N,x}$ is given by
\begin{equation}\label{4croMNx}
\langle M^{N,x}\rangle_t=\int_0^t \left\lbrace 4Z^{N,x}_r+\frac{1}{N}\vert\vert f\vert\vert_{N,0,Z^{N,x}_r}\right\rbrace dr,
\end{equation}
where for any $z=\frac{k}{N}$, $z'=\frac{k'}{N}$,  $k\in\Z_+$ such that \ $k\le k'$,
  $$\vert\vert f\vert\vert_{N,z,z'}=\sum_{i=k+1}^{k'}\vert f(\frac{i}{N})-f(\frac{i-1}{N})\vert .$$
Now we precise the law of the pair $\left(Z^{N,x},Z^{N,y}\right) $, for any $x,y\in\mathbb{R}_+$ such that $x\le y$. Let $V^{N,x,y}:=Z^{N,y}-Z^{N,x}$,  and consider the pair of process $\left(Z^{N,x},V^{N,x,y}\right) $, which starts from $\left(\frac{\lfloor Nx\rfloor}{N} ,\frac{\lfloor Ny\rfloor-\lfloor Nx \rfloor}{N} \right)$ at time $t=0$, and whose dynamic is described by: 
$\left( Z^{N,x},V^{N,x,y}\right) $ jumps  \\
 from $(\frac{i}{N},\frac{j}{N})$ to
$\begin{cases} (\frac{i+1}{N},\frac{j}{N})\text{ at rate }  2Ni +\sum_{k=1}^{i}(f(\frac{k}{N})-f(\frac{k-1}{N})^+\\ 
(\frac{i-1}{N},\frac{j}{N}) \text{ at rate } 2Ni+\sum_{k=1}^{i}(f(\frac{k}{N})-f(\frac{k-1}{N}))^-\\
(\frac{i}{N},\frac{j+1}{N}) \text{ at rate } 2Nj +\sum_{k=1}^{j}(f(\frac{i+k}{N})-f(\frac{i+k-1}{N})^+\\
(\frac{i}{N},\frac{j-1}{N}) \text{ at rate } 2Nj+\sum_{k=1}^{j}(f(\frac{i+k}{N})-f(\frac{i+k-1}{N})^-
\end{cases}$.\\
The process $V^{N,x,y}$ can be expressed as follows. 
\begin{align}\label{4fVN}
V^{N,x,y}_t=\frac{\lfloor Ny\rfloor -\lfloor Nx\rfloor}{N}+\int_0^t \left[f(Z^{N,x}_r+V^{N,x,y}_r)-f(Z^{N,x}_r)\right] dr+M^{N,x,y}_t,
\end{align}
where $M^{N,x,y}$ is a local martingale whose conditional quadratic variation  $\langle M^{N,x,y}\rangle$ is given by
\begin{align}\label{4croMNxy}
\langle M^{N,x,y}\rangle_t =\int_0^t \left\lbrace 4V^{N,x,y}_r+\frac{1}{N}\vert\vert f\vert\vert_{N,Z^{N,x}_r,V^{N,x,y}+Z^{N,x}_r}\right\rbrace dr.
\end{align}
Since $Z^{N,x}$ and $V^{N,x,y}$ never jump at the same time, 
\begin{align} 
[M^{N,x}, M^{N,x,y}]=0, \text{ hence } \langle M^{N,x},M^{N,x,y}\rangle =0,\label{4croch}
\end{align}
which implies that the martingales $M^{N,x}$ and $M^{N,x,y}$ are orthogonal.

Consequently, $Z^{N,x}+V^{N,x,y}$ solves the SDE
\begin{align*}
Z^{N,x}_t+V^{N,x,y}_t=&\frac{\lfloor Ny\rfloor}{N}+\int_0^t f(Z^{N,x}_r+V^{N,x,y}_r)dr+\tilde{M}^{N,x,y}_t.
\end{align*}
where $\tilde{M}^{N,x,y}$ is a local martingale with $\langle \tilde{M}^{N,x,y}\rangle$ given by
\begin{align*}
\langle \tilde{M}^{N,x,y}\rangle_t &=\langle M^{N,x}\rangle_t+\langle M^{N,x,y}\rangle_t = \langle M^{N,x+y}\rangle_t,\quad \forall t\ge 0.
\end{align*}
We then deduce that for any $x,y\in\R_+$ such $x\le y, $
\[Z^{N,x}+V^{N,x,y}\stackrel{(d)}{=}Z^{N,y}.\] 
%It is not hard to see that in fact conditionally upon $\left{Z^{N,x'},x'\le x\right}$, $M^{N,x,y}$ is a local martingale.
In fact, we have that 
\begin{align*}
V^{N,x,y}_t&=\frac{\lfloor Ny\rfloor -\lfloor Nx\rfloor}{N}+\frac{1}{N}P^1\left(N\int_0^t\sum_{k=1}^{NV^{N,x,y}_s}\left(f(Z^{N,x}_s+\frac{k}{N})-f(Z^{N,x}_s+\frac{k}{N})\right)^+ds\right)\\
&-\frac{1}{N}P^2\left(N\int_0^t\sum_{k=1}^{NV^{N,x,y}_s}\left(f(Z^{N,x}_s+\frac{k}{N})-f(Z^{N,x}_s+\frac{k}{N})\right)^-ds\right)\\
&+\frac{1}{N}P^3\left(2N^2\int_0^tV^{N,x,y}_sds\right)-\frac{1}{N}P^4\left(2N^2\int_0^tV^{N,x,y}_sds\right),
\end{align*}
where $P^1$, $P^2$, $P^3$ and $P^4$ are mutually independent standard Poisson processes which are all independent of $Z^{N,x}_.$. 
It follows that conditionally upon  $\left\lbrace Z^{N,x'},x'\le x\right\rbrace$, $M^{N,x,y}$ is a local martingale.
\begin{remark}
We can also proceed to the renormalisation of the exploration process to provide a discrete Ray Knight representation of the process $Z^{N,x}$.
We choose the slope $p=2N$ and we denote by $H^N$ the exploration processes associated to the forest $\mathcal{F}^{N,x}$ of $\lfloor Nx\rfloor$ trees. Let  $L^N_s(t)$  be the local time of the process  $H^N$ at level $t$ up to time $s$. 
At any time $s$, the rate of minima of $H^N$  is equal to
$$4N^2+N\left[f(\frac{\lfloor NL^N_s(H^N_s)\rfloor}{N}+1/N)-f(\frac{\lfloor NL^N_s(H^N_s)\rfloor}{N})\right]^+,$$
and the rate of maxima is equal to $$4N^2 +N\left[f(\frac{\lfloor NL^N_s(H^N_s)\rfloor}{N}+1/N)-f(\frac{ \lfloor NL^N_s(H^N_s)\rfloor }{N})\right]^-.$$ Let  $S^{N,x}$ be the time to explore the forest $\mathcal{F}^{N,x}$. We have that 
$$S^{N,x}=\inf\lbrace s>0; L^N_s(0)\ge \frac{\lfloor Nx\rfloor}{N}\rbrace .$$
Under the assumption that $S^{N,x}<\infty$ $a.s.$ for all $x>0$, the discrete Ray Knight representation with the renormalization becomes: 
\begin{align*}
\left( Z^{N,x}_t, t\ge 0, x\ge 0\right)\equiv\left( L^N_{S^{N,x}}(t),t\ge 0, x\ge 0\right).
\end{align*}
One could probably deduce from this  discrete approximation the Ray Knight representation of the general Feller diffusion by a limiting argument, as it is done in \cite{4BPS} in the linear case and in \cite{4LPW} in the quadratic case.  But in this work  we use stochastic analysis tools for proving our extended Ray Knight theorem. 
\end{remark}
\subsection{Continous model with a general competition}\label{4contm}
Given a space-time white noise $W(ds,du)$, we now define an $\R_+$--valued two--parameter stochastic process 
$\lbrace Z^{x}_t,\ t\ge 0,$ $x\ge 0\rbrace$ which is such that for each fixed $x>0$, $\left\lbrace Z^{x}_t,t\ge 0\right\rbrace $  is a continuous process, solution of the SDE \eqref{4fA}. 
We have that for any $0<x<y$, $\{V^{x,y}_t:=Z^{y}_t-Z^x_t,\ t\ge0\}$
solves the SDE
 \begin{align} \label{4fAV}
V^{x,y}_t= y-x+\int_0^t \left[ f(Z^x_s+V^{x,y}_s) -f(Z^x_s)\right] ds+2\int_0^t\int_{Z^x_s}^{Z^x_s+V^{x,y}_s}W(ds,du)
\end{align}
The process $V^{x,y}$ is nonnegative  almost surely. We have that $\int_0^t\int_{0}^{Z^x_s}W(ds,du)$ and $\int_0^t\int_{Z^x_s}^{Z^x_s+V^{x,y}_s}W(ds,du)$ are orthogonal since $$[0,Z^x_s]\cap(Z^x_s, Z^x_s+V^{x,y}_s]=\emptyset, $$
and 
$$\int_0^t\int_{0}^{Z^x_s}W(ds,du)+\int_0^t\int_{Z^x_s}^{Z^x_s+V^{x,y}_s}W(ds,du)=\int_0^t\int_{0}^{Z^x_s+V^{x,y}_s}W(ds,du)\text{ a.s.}$$
This implies that  $Z^{y}=Z^x+V^{x,y}$ a.s.  It follows that,  for each $t\ge 0$,  the process $\{Z^x_t, x\ge0\}$ is almost surely non decreasing and  for $0\le x<y$, the conditional law of $Z^y_\cdot$, given $\{Z^{x'}_t, x'\le x, t\ge 0\}$ and  $Z^x_t=z(t)$, $t\ge0$, is the law of the sum of $z$ plus the solution of \eqref{4fAV} with $Z^x_t$ replaced by $z(t)$. Note that when $Z^x_.$ is replaced by a deterministic trajectory $z$, the solution of \eqref{4fAV} is independent of $\{Z^{x'}_., x'< x \}$ .  Hence the process $\{Z^x_\cdot,\ x\ge0\}$ is a Markov process with values in $C([0,\infty),\R_+)$, the space of continuous functions from $[0,\infty)$ into $\R_+$, starting from $0$ at $x=0$. In the case $f$ linear, the increments of the mapping $x\to Z^x_t$ are independent, for each $t>0$. 

For $x\ge 0$, define $T^x_0$ the extinction time of the process $Z^x$ by:
$$T^x_0=\inf \left\lbrace t>0; Z^x_t=0\right\rbrace. $$ 
For any $x\ge 0$, we call the process $Z^x$ subcritical if it goes extinct almost surely in finite time i.e if $T^x_0$ is finite almost surely. The assumption \textbf{A} implies that $\frac{f(x)}{x}$ is bounded. Let us introduce the notation
\begin{equation}\label{lambda}
\Lambda(f):=\int_1^{\infty}\exp\left(-\frac{1}{2}\int_{1}^u \frac{f(r)}{r}dr\right)du.
\end{equation}
We have the following Proposition.
\begin{proposition}\label{subcritic}
Suppose that  $f$ satisfies hypothesis \textbf{A}. For any $x\ge 0$, $Z^x$ is subcritical if and only if 
$\Lambda(f)=\infty$.
In particular we have:
\begin{itemize} 
\item[i)]A sufficient condition for $\P\left(T^x_0<\infty\right)=1$ is: there exists $z_0\ge 1$  such that $f(z)\le 2$, $\forall z\ge z_0$, 
\item[ii)] A sufficient condition for $\P\left(T^x_0=\infty\right)>0$ is: there exists $z_0>1$ and $\delta>0$ such that  $f(z)\ge 2+\delta$, $\forall$ $z\ge z_0$. 
\end{itemize}
\end{proposition}
\bpf 

%For any $y\ge 0$, we define more generally $$T^x_y=\inf \left\lbrace t>0; Z^x_t=y \right\rbrace .$$ 
%First of all, for any $0\le x\le x_0$, $\P \left(T^x_0\wedge T^x_{x_0+1}<\infty\right)=1$.
%Consequently if $Z^x$ never reaches the level $x_0+1$, $T^x_0<\infty$ $a.s.$. Next, it is plain that $$\inf_{0\le x\le x_0}  \P \left(T^x_0<T^x_{x_0+1}\right):=p>0 .$$
%Let $Z^x_t$ start afresh from $x_0+1$. On the interval [0, $T^{x_0+1}_{x_0}$], $Z^{x_0+1}$ is bounded from above by the martingale $$2\int_0^t\sqrt{Z^{x_0+1}_r}dWr$$
%which is a time  changed Brownian motion, hence it reaches $x_0$ in finite time. Consequently $T^{x_0+1}_{x_0}<\infty $ $a.s.$ . Each time $Z^x_t$ reaches $x_0$, from the strong Markov property it has a probability $p$ of reaching $0$ before the level $x_0+1$, independently of the past. The result is established. 
Let $S\in C^2(\R_+)$ and $0\le a< x< b$. By It\^o's formula applied to the process $Z^x$ and the function $S$, we have that for any  $t\ge0$,
\begin{align}\label{itoz}
S(Z^x_t)=S(x)+\int_0^t\left(S'(Z^x_s)f(Z^x_s)+2S''(Z^x_s)Z^x_s\right)ds+2\int_0^tS'(Z^x_s)\sqrt{Z^x_s}dW_s.
\end{align}
Let us denote by $\mathcal{A}$ the generator of $Z^x$. If we can find a strictly increasing function $S$ on the interval $[a,b]$ such that $\mathcal{A} S\equiv 0$, then the drift term in \eqref{itoz}  vanishes and so $Z^x$ will be just a time changed Brownian motion in [$S(a)$, $S(b)$]. Such a function $S$ is called a scale function of the diffusion $Z^x$. We choose as scale function:  for any $z\ge 0$,
$$S(z)=\int_1^z\exp\left(-\frac{1}{2}\int_1^u \frac{f(r)}{r}dr\right)du.$$
Let us denote by $T^x_y$ the random time at which $Z^x$ hits $y$ for the first time. We have for any $0\le a< x< b$\begin{align*}
\P(T^x_a<T^x_b)=\frac{S(b)-S(x)}{S(b)-S(a)}, \text{ and } \P(T^x_a<\infty)=\lim_{b\to \infty}\P(T^x_a<T^x_b).
\end{align*}
If the function $S(z)$ tends to infinity as $z$ goes to infinity, then $\P(T^x_a<\infty)=1$. Otherwise $0<\P(T^x_a<\infty)<1$.
From this we deduce that $Z^x$ goes extinct almost surely in finite time if and only if $\lim_{z\to \infty}S(z)=\infty$ i.e. if and only if 
$\Lambda(f) =\infty$.
The rest of the Proposition is immediate.
%Now suppose that there exists $z_0\ge 0$ such that $f(z)\le 2$, $\forall$ $z\ge z_0$. \\
%We have that
%\begin{align*}
%\int_0^{\infty}\exp\left(-\frac{1}{2}\int_{z_0}^u \frac{f(r)}{r}dr\right)du  \ge \int_0^{\infty}\frac{1}{u}du =\infty.
%\end{align*}
%Hence a sufficient condition for $T^x_0$ to be finite almost surely is: there exists $z_0\ge 0$ such that $f(z)\le 2$, $\forall$ $z\ge z_0$.  

%Let suppose now that  there exists $z_0\ge 0$ and $\delta>0$ such that $f(z)\ge 2+\delta$, $\forall$ $z\ge z_0$. We then have
%\begin{align*}
%\int_0^{\infty}\exp\left(-\frac{1}{2}\int_{z_0}^u \frac{f(r)}{r}dr\right)du & \le \int_0^{\infty}\frac{1}{u^{1+\delta}}du \\
%& < \infty.
%\end{align*}
%A sufficient condition for $\P\left( T^x_0=\infty \right)>0$ is: there exists $z_0\ge 0$ and $\delta>0$ such that $f(z)\ge 2+\delta$, $\forall$ $z\ge z_0$.
\epf
\section{Convergence as $N\to\infty$}
The aim of this section is to prove the convergence in law as $N\to\infty$ of the two--parameter process 
$\{Z^{N,x}_t,\ t\ge0, x\ge0\}$ defined in section \ref{4secRenorm} towards the process $\{Z^{x}_t,\ t\ge0, x\ge0\}$ defined in section \ref{4contm}. We need to make precise the topology for which this convergence will hold.
We note that the process $Z^{N,x}_t$ (resp. $Z^x_t$) is a Markov processes indexed by $x$,
with values in the space of c\`adl\`ag (resp. continuous) functions of $t$
$D(([0,\infty);\R_+)$ (resp. $C(([0,\infty);\R_+)$). So it will be natural to consider a topology of functions of $x$, with values in functions of $t$.

For each fixed $x$, the process $t\to Z^{N,x}_t$ is c\`adl\`ag, constant between its jumps, with jumps of size $\pm N^{-1}$, while the limit process $t\to Z^{x}_t$ is continuous.
On the other hand, both $Z^{N,x}_t$ and $Z^{x}_t$ are discontinuous as functions of $x$. $x\to Z^x_\cdot$ has countably many jumps on any compact interval, but the mapping $x\to\{Z^x_t,\ t\ge\epsilon\}$, where
$\epsilon>0$ is arbitrary,  has finitely many jumps on any compact interval, and it is constant between its jumps. 
%(in this statement we exclude the cases ($\alpha<1, \theta>0$
%and $\alpha=1, \theta>\gamma$)
Recall that $D([0,\infty);\R_+)$ equipped with the distance $d^0_{\infty}$ defined by $(16.4)$ in \cite{4BI} is separable and complete, see Theorem $16.3$ in \cite{4BI}. We have the following statement
\begin{theorem}\label{4thConv} Suppose that the Hypothesis $\bf{A}$ is satisfied. Then as $N\to\infty$,
$$\{Z^{N,x}_t,\ t\ge0, x\ge0\}\Rightarrow\{Z^{x}_t,\ t\ge0, x\ge0\}$$
in $D([0,\infty);D([0,\infty);\R_+))$, equipped with the Skohorod topology of the space of 
c\`adl\`ag functions of $x$, with values in the Polish space $D([0,\infty);\R_+)$ equipped with the metric $d^0_{\infty}$.
\end{theorem}
\subsection*{Proof of the theorem}
To prove the theorem, we first show that for fixed $x\ge 0$ the sequence $\left\lbrace Z^{N,x},~N\geq 0\right\rbrace$ is tight in $D([0,\infty);\R_+)$.
\subsection{Tightness of $Z^{N,x}$}
For this end, we first establish a few lemmas.
\begin{lemma} \label{4lm1}
For all $T>0$, $x\ge 0$, there exist a constant $C_0>0$ such that for all $N\ge 1$,
\begin{align*} \sup_{0\leq t\leq T}\E\left( Z^{N,x}_t\right) \le C_0.
\end{align*}
Moreover, for all $t\ge 0$, $N\ge1$,
\begin{align*}
\E \left(-\int_0^t f(Z^{N,x}_r)dr\right) \le x.
\end{align*}
\end{lemma}
\bpf
Let $\left(\tau_n,n\ge 0\right)$ be a sequence of stopping times such that $\tau_n$ tends to infinity as $n$ goes to infinity and  for any $n$, $\left(M^{N,x}_{t\wedge \tau_n},t\ge 0\right)$ is a martingale and $Z^{N,x}_{t\wedge\tau_n}\le n$. Taking the expectation on both sides  of equation \eqref{4fAN} at time $t\wedge\tau_n$, we obtain
\begin{align}\label{4equ1}
\E\left(Z^{N,x}_{t\wedge \tau_n}\right)=\frac{\lfloor Nx\rfloor}{N}+\E\left(\int_0^{t\wedge \tau_n} f(Z^{N,x}_r)dr\right).
\end{align}
It follows from the hypothesis \textbf{A} on $f$ that
\begin{align*}
\E\left(Z^{N,x}_{t\wedge \tau_n}\right) \le \frac{\lfloor Nx\rfloor}{N}+\beta \int_0^t\E (Z^{N,x}_{r\wedge \tau_n})dr
\end{align*}
From Gronwall and Fatou Lemmas, we deduce that there exists a constant  $C_0>0$ which depends only upon $x$ and $T$ such that
\[\sup_{N\ge 1}\sup_{0\le t\le T}\E\left(Z^{N,x}_t\right)\le C_0.\]
From \eqref{4equ1}, we deduce that
$$-\E\left(\int_0^{t\wedge \tau_n} f(Z^{N,x}_r)dr\right)\le \frac{\lfloor Nx\rfloor}{N}  .$$
Since $-f(Z^{N,x}_r )\ge -\beta Z^{N,x}_r$ , the second statement follows using Fatou's Lemma and the first statement.
\epf

We now have the following lemma.
\begin{lemma}\label{4lemcro}
For all $T>0$, $x\ge 0$, there exists a constant $C_1>0$ such that : 
$$\sup_{N\geq 1}\E\left(\langle M^{N,x}\rangle_T\right)\le C_1.$$
\end{lemma}
\bpf
For any $N\ge 1$ and $k,k' \in\Z_+$ such that $k\le k'$, we set $z=\frac{k}{N}$ and $z'=\frac{k'}{N}$. We deduce from hypothesis \textbf{A} on $f$ that 
\begin{align*}
\vert \vert f\vert\vert_{N,z,z'}&=\sum_{i=k+1}^{k'}\left\lbrace\left( f(\frac{i}{N})-f(\frac{i-1}{N})\right)^++\left( f(\frac{i}{N})-f(\frac{i-1}{N})\right)^-\right\rbrace\\
&=\sum_{i=k+1}^{k'}\left\lbrace 2\left( f(\frac{i}{N})-f(\frac{i-1}{N})\right)^+-\left( f(\frac{i}{N})-f(\frac{i-1}{N})\right)\right\rbrace.
\end{align*}
Consequentely,
\begin{align}\label{4bornvaria}
\vert \vert f\vert\vert_{N,z,z'} \le 2\beta(z'-z)+f(z)-f(z').
\end{align}
We deduce from \eqref{4bornvaria}, \eqref{4croMNx} and Lemma \ref{4lm1} that 
\begin{align*}
\E\left(\langle M^{N,x}\rangle_T\right)& \le \int_0^T \left\lbrace (4+\frac{2\beta}{N})\E (Z^{N,x}_r)- \frac{1}{N} \E\left(f(Z^{N,x}_r\right)\right\rbrace dr\\
&\le (4+\frac{2\beta}{N})C_0T+\frac{x}{N}.
\end{align*}
Hence the lemma.
\epf
It follows from this that $M^{N,x}$ is in fact a square integrable martingale. We also have 
\begin{lemma}\label{4lm2}
For all $T>0$, $x\ge 0$, there exist two  constants $C_2,C_3>0$ such that ~: 
\begin{align*}  \sup_{N\geq 1}\sup_{0\leq t\leq T}\E\left(  Z^{N,x}_t\right)^2 &\le C_2,\\
\sup_{N\geq 1}\sup_{0\leq t\leq T}\E\left(  -\int_0^t Z^{N,x}_r f(Z^{N,x}_r)dr \right)& \leq C_3.
\end{align*}
\end{lemma}
\bpf
We deduce from \eqref{4fAN} and It\^o's formula that
\begin{align}\label{4carre}
\left(Z^{N,x}_t\right)^2 = \left(\frac{ \lfloor Nx\rfloor }{N}\right)^2+2\int_0^t Z^{N,x}_r f(Z^{N,x}_r)dr+\langle M^{N,x}\rangle_t +M^{N,x,(2)}_t,
%\nonumber \\
%&+ \int_0^t\left[4Z^{N,x}_r+\frac{1}{N}\vert \vert f \vert\vert_{N,0,Z^{N,x}_r}\right]dr+M^{N,x,(2)}_t,
\end{align}
where $M^{N,x,(2)}$ is a local martingale. 
Let $\left( \sigma_n,n\ge 1\right)$ be a sequence of stopping times such that $\lim_{n\to\infty}\sigma_n=+\infty$ and for each $n\ge 1$, $\left( M^{N,x,(2)}_{t\wedge\sigma_n}, t\ge 0\right)$ is a martingale. Taking the expectation on the both sides of \eqref{4carre} at time $t\wedge \sigma_n$ and using  hypothesis \textbf{A}, Lemma \ref{4lemcro}, Gronwall and  Fatou lemmas  we  obtain that for all $T>0$, there exists a constant $C_2>0$ such that : 
$$\sup_{N\geq 1}\sup_{0\leq t\leq T}\E\left(  Z^{N,x}_t\right)^2 dr \leq C_2.$$
We also have that
\begin{align*}
2\E\left(-\int_0^{t\wedge\sigma_n} Z^{N,x}_r f(Z^{N,x}_r)dr \right) \le  \left(\frac{ \lfloor Nx\rfloor }{N}\right)^2+ C_1
\end{align*}
From Hypothesis $\textbf{A}$, we have $- Z^{N,x}_r f(Z^{N,x}_r)\ge - \beta(Z^{N,x}_r)^2$. The result now follows from Fatou's Lemma.
\epf 

We want to check tightness of the sequence $\left\lbrace Z^{N,x},~N\geq 0\right\rbrace $ using Aldous' criterion.  Let $\left\lbrace \tau_N,~N\geq 1\right\rbrace$ be a sequence of stopping time in $[0,T]$. We deduce from Lemma \ref{4lm2}
\begin{proposition} \label{4prop2}
For any $T>0$ and $\eta$, $\epsilon >0$, there exists $\delta>0$ such that 
\[\sup_{N\ge 1}\sup_{0\le \theta\le \delta}\P\left(\left\vert  \int_{\tau_N}^{(\tau_N+\theta)\wedge T} f(Z^{N,x}_r)dr\right\vert \ge \eta \right)\le \epsilon . \]
\end{proposition}
\bpf
Let $c$ be a non negative constant. We have
\begin{align*}
\left\vert \int_{\tau_N}^{(\tau_N+\theta)\wedge T}f(Z^{N,x}_r)dr\right\vert &\le \sup_{0\le r \le c }\vert f(r)\vert\delta +\int_{\tau_N}^{\tau_N+\theta}{\bf1}_{ \lbrace Z^{N,x}_r>c\rbrace }\vert f(Z^{N,x}_r)\vert dr
\end{align*}
But 
\begin{align*}
\int_{\tau_N}^{\tau_N+\theta}{\bf1}_{ \lbrace Z^{N,x}_r>c\rbrace }\vert f(Z^{N,x}_r)\vert dr & \le c^{-1}\int_0^T Z^{N,x}_r \left(f^+(Z^{N,x}_r)+f^-(Z^{N,x}_r)\right)dr\\
&\le c^{-1}\int_0^T\left(2Z^{N,x}_r f^+(Z^{N,x}_r)-Z^{N,x}_rf(Z^{N,x}_r) \right)dr \\
&\le c^{-1}\int_0^T\left(2\beta (Z^{N,x}_r)^2-Z^{N,x}_rf(Z^{N,x}_r) \right)dr.
\end{align*} 
From this and Lemma \ref{4lm2}, we deduce  that $\forall$ $N\ge 1$
\begin{align*}
\sup_{0\le \theta\le \delta}\P\left(\Big\vert  \int_{\tau_N}^{(\tau_N+\theta)\wedge T} f(Z^{N,x}_r)dr\Big\vert\ge \eta \right) &\le \eta^{-1}\E\left(\left\vert  \int_{\tau_N}^{(\tau_N+\theta)\wedge T}
f(Z^{N,x}_r)dr\right\vert \right)\\
& \le \sup_{0\le r \le c}\frac{\vert f(z)\vert\delta}{ \eta }+\frac{A}{c\eta},
\end{align*}
with $A=2\beta C_2T+C_3$. The result follows by choosing $c=2A/\epsilon\eta$, and then $\delta=\epsilon\eta/2\sup_{0\le r\le c}\vert f(z)\vert $.
\epf
From Proposition \ref{4prop2}, the Lebesgue integral term in the right hand side of \eqref{4fAN} satisfies Aldou's condition $[A]$, see \cite{AD}. The same Proposition, Lemma \ref{4lm1}, \eqref{4croMNx} and \eqref{4bornvaria} imply that $<M^{N,x}>$ satisfies the same condition, hence so does $M^{N,x}$, according to Rebolledo's theorem, see \cite{4JM}. Since all jumps are of size $\frac{1}{N}$, tightness follows. We have proved
%From  \eqref{4fAN}, Proposition \ref{4prop2}, Lemma \ref{4lm1}, \eqref{4croMNx} and \eqref{4bornvaria} we deduce that
\begin{proposition}\label{4tightxfixed}
For any fixed $x\ge 0$, the sequence of processes $\left\lbrace Z^{N,x},~N\ge 1\right\rbrace $ is tight in $D\left([0,\infty);\R_+ \right) $.
\end{proposition}
We deduce from Proposition \ref{4tightxfixed} the following Corollary.
\begin{corollary}\label{4cortight}
For any $0\le x<y$ the sequence of processes $\left\lbrace V^{N,x,y},~N\ge 1\right\rbrace$ is tight in $D\left([0,\infty);\R_+ \right) $
\end{corollary}
\bpf 
For  any $x$ fixed  the process $Z^{N,x}$ has jumps equal to $\pm\frac{1}{N}$ which tends to zero as $N\to\infty$. It follows from that and  equation \eqref{4fAN} that any weak limit of a converging subsequence of $Z^{N,x}$ is continuous and is the unique weak solution of equation \eqref{4fA}. We deduce that for any $x,y\ge 0$, 
the sequence  $\left\lbrace Z^{N,y}-Z^{N,x},N\ge 1\right\rbrace$ is tight since $\left\lbrace Z^{N,x},N\ge 1\right\rbrace$ and$\left\lbrace Z^{N,y},N\ge 1\right\rbrace$ are tight and both have a continuous limit as $N\to\infty$.
\epf 
\subsection{Proof of Theorem \ref{4thConv}}
From Theorem 13.5 in \cite{4BI}, Theorem \ref{4thConv} follows from the two next Propositions
 
\begin{proposition}\label{4findimconv}
For any $n\in\N$, $0\le x_1<x_2<\cdots <x_n$,  $$\left( Z^{N,x_1},Z^{N,x_2},\cdots,Z^{N,x_n}\right)\Rightarrow \left( Z^{x_1},Z^{x_2},\cdots,Z^{x_n}\right)$$
as $N\to\infty$, for the topology of locally uniform convergence in $t$.
\end{proposition}
\bpf
We prove the statement in the case $n=2$ only. The general statement can be proved in a very similar way.
For $0\le x_1<x_2$, we consider the process $\left(Z^{N,x_1},V^{N,x_1,x_2}\right)$, using the notations from section \ref{4discrete}.   The argument preceding the statement of Proposition \ref{4tightxfixed} implies that  the sequences of martingales $M^{N,x_1}$ and $M^{N,x_1,x_2}$ are tight. Hence  \\
$\left( Z^{N,x_1},V^{N,x_1,x_2}, M^{N,x_1},M^{N,x_1,x_2}\right)$ is tight. Thanks  to \eqref{4fAN}, \eqref{4fVN}, \eqref{4croMNx}, \eqref{4croMNxy} and \eqref{4croch}, any converging subsequence of  \\
$\left\lbrace Z^{N,x_1},V^{N,x_1,x_2}, M^{N,x_1},M^{N,x_1,x_2}, N\ge 1\right\rbrace $ has a weak limit \\
 $\left(Z^{x_1},V^{x_1,x_2}, M^{x_1},M^{x_1,x_2}\right)$  which satisfies 
 \begin{align*}
 Z^{x_1}_t&=x_1+\int_0^tf(Z^{x_1}_s)ds+M^{x_1}_t\\
 V^{x_1,x_2}_t & =x_2-x_1+\int_0^tf\left[f(Z^{x_1}_s+V^{x_1,x_2}_s)-f(Z^{x_1}_s)\right]ds+M^{x_1,x_2}_t,
 \end{align*}
 where the continuous martingales $M^{x_1}$ and $M^{x_1,x_2}$  satisfy 
\begin{align*}
\langle M^x\rangle_t=4\int_0^t Z^{x_1}_sds, ~\langle M^{x_1,x_2}\rangle_t=4\int_0^tV^{x_1,x_2}_sds, ~
\langle M^{x_1},M^{x_1,x_2}\rangle_t =0.
\end{align*} 
This implies that the pair $\left(Z^{x_1},V^{x_1,x_2}\right)$ is a weak solution of the system of SDEs \eqref{4fA} and \eqref{4fAV}, driven by the same space-time white noise. The result follows from the uniqueness of the system, see \cite{4DL}.
\epf

\begin{proposition}\label{4tightD}
There exists a constant $C$, which depends only upon $\theta$ and $T$, such that
for any $0\le x<y<z$, which are such that $y-x\le1$, $z-y\le 1$, 
$$\E\left[\sup_{0\le t\le T}\vert Z^{N,y}_t-Z^{N,x}_t\vert^2\times
\sup_{0\le t\le T}\vert Z^{N,z}_t-Z^{N,y}_t\vert^2\right]\le C|z-x|^2.$$
\end{proposition}
%Now we let us prove the tightness of the sequence $\left\lbrace Z^{N,x},x\ge 0\right\rbrace _{N\ge 1}$ %in $D([0,\infty);C([0,T]))$ for all $T>0$. For this end we show the following inequality which yields the %tightness. 

We first prove the 
\begin{lemma} \label{4lmVN} For any $0\le x<y$, we have
$$\sup_{0\le t\le T}\E \left(Z^{N,y}_t-Z^{N,x}_t\right)=\sup_{0\le t\le T}\E(V^{N,x,y}_t)\le \left(\frac{\lfloor Ny \rfloor}{N}-\frac{\lfloor Nx \rfloor }{N}\right)e^{\beta T},$$
\end{lemma}
\bpf
Let $\left( \sigma_n,n\ge 0\right)$ be a sequence of stopping times such that $\lim_{n\to\infty}\sigma_n=+\infty$ and $\left( M^{N,x,y}_{t\wedge\sigma_n}\right)$ is a martingale. Taking the expectation on the both sides of \eqref{4fVN} at time $t\wedge \sigma_n$ we obtain that
\begin{align}\label{4vb}
\E(V^{N,x,y}_{t\wedge \sigma_n})\le\left(\frac{\lfloor Ny \rfloor}{N}-\frac{\lfloor Nx \rfloor }{N}\right)+\beta\int_0^t \E (V^{N,x,y}_{r\wedge \sigma_n})dr 
\end{align}
Using Gronwall and Fatou lemmas, we obtain that $$\sup_{0\le t\le T}\E(V^{N,x,y}_t)\le\left(\frac{\lfloor Ny \rfloor}{N}-\frac{\lfloor Nx \rfloor }{N}\right)e^{\beta T}.$$ 
\epf
\noindent{\sc Proof of Proposition \ref{4tightD}}
Using equation \eqref{4fVN}, a stopping time argument as above, Lemma \ref{4lmVN} and  Fatou's lemma, where we take advantage of the inequality  $f(Z^{N,x}_r)-f(Z^{N,x}_r+V^{N,x,y}_r)\ge -\beta V^{N,x,y}_r$, we deduce that
\begin{align}\label{4ine2}
\E\left(\int_0^t \left[f(Z^{N,x}_r)-f(Z^{N,x}_r+V^{N,x,y}_r) \right]dr\right) \le \frac{\lfloor Ny \rfloor}{N}-\frac{\lfloor Ny \rfloor }{N}.
\end{align}
We now deduce from \eqref{4croMNxy}, Lemma \ref{4lmVN}, inequalities \eqref{4ine2} and \eqref{4bornvaria} that for each $t>0$, there exists a constant $C(t)>0$ such that
\begin{align}\label{4corbor}
\E\left( \langle M^{N,x,y}\rangle_t \right) \le C(t)\left(\frac{\lfloor Ny \rfloor}{N}-\frac{\lfloor Ny \rfloor }{N}\right). 
\end{align}
This implies that $M^{N,x,y}$ is in fact a square integrable  martingale.
For any $0\le x<y<z$, we have $Z^{N,z}_t-Z^{N,y}_t=V^{N,y,z}_t$ and $Z^{N,y}_t-Z^{N,x}_t=V^{N,x,y}_t$ for any $t\ge0$. On the other hand we deduce from \eqref{4fVN} and the hypothesis $\bf{A}$
\begin{align*} 
\sup_{0\le t\le T}(V^{N,x,y}_t )^2 & \le  3\left(\frac{\lfloor Ny \rfloor}{N}-\frac{\lfloor Nx \rfloor }{N}\right)^2+3\beta^2T\int_0^T \sup_{0\le s\le r}(V^{N,x,y}_s)^2dr \\
&+3\sup_{0\le t\le T}\left(M^{N,x,y}_t\right)^2,
\end{align*}
\begin{align*} 
\sup_{0\le t\le T}(V^{N,y,z}_t )^2 & \le 3\left(\frac{\lfloor Nz \rfloor}{N}-\frac{\lfloor Ny \rfloor }{N}\right)^2+3\beta^2T\int_0^t \sup_{0\le s\le r}(V^{N,y,z}_s)^2dr \\
&+3\sup_{0\le t\le T}\left(M^{N,y,z}_t\right)^2.
\end{align*}
Now let $\mathcal{G}^{x,y}:=\sigma\left( Z^{N,x}_t, Z^{N,y}_t,t\ge 0\right)$ be the filtration generated by $Z^{N,x}$ and $Z^{N,y}$. It is clear that  for any $t$, $V^{N,x,y}_t$ is measurable with respect to $\mathcal{G}^{x,y}$. We then have 
\begin{align*}
\E\left[\sup_{0\le t\le T}\vert V^{N,x,y}_t\vert^2\times \sup_{0\le t\le T}\vert V^{N,y,z}_t\vert^2\right]=\E\left[\sup_{0\le t\le T}\vert V^{N,x,y}_t\vert^2 \E\left(\sup_{0\le t\le T}\vert V^{N,y,z}_t\vert^2 \vert \mathcal{G}^{x,y}\right) \right].
\end{align*}
Conditionally upon $Z^{N,x}$ and $Z^{N,y}=u(.)$, $V^{N,y,z}$ solves the following SDE
\begin{align*}
V^{N,y,z}_t=\frac{\lfloor Nz\rfloor -\lfloor Ny\rfloor}{N}+\int_0^t \left[f(V^{N,y,z}_r+u(r))-f(u(r))\right]dr+M^{N,y,z}_t,
\end{align*}
where $M^{N,y,z}$ is a martingale conditionally upon $\mathcal{G}^{x,y}$, hence the arguments used in Lemma \ref{4lmVN} lead to
$$\sup_{0\le t\le T}\E\left(V^{N,y,z}_t\vert \mathcal{G}^{x,y}\right) \le\left(\frac{\lfloor Nz \rfloor}{N}-\frac{\lfloor Ny \rfloor }{N}\right)e^{\beta T}, $$ 
and those used to prove \eqref{4ine2} yield
$$\E\left(\int_0^t f(Z^{N,y}_r)-f(Z^{N,y}_r+V^{N,y,z}_r) dr\vert  \mathcal{G}^{x,y}\right) \le \frac{\lfloor Nz \rfloor}{N}-\frac{\lfloor Ny \rfloor }{N}.$$
From this we deduce (see the proof of \eqref{4corbor}) that
$$\E\left( \langle M^{N,y,z}\rangle_t \vert \mathcal{G}^{x,y}\right) \le C(t)\left(\frac{\lfloor Nz \rfloor}{N}-\frac{\lfloor Ny \rfloor }{N}\right).$$
From Doobs's inequality we have
\begin{align*}
\E\left(\sup_{0\le t\le T}\vert M^{N,y,z}_t\vert^2 \vert \mathcal{G}^{x,y}\right)&=\E\left( \langle M^{N,y,z}\rangle_T\vert \mathcal{G}^{x,y}\right) \\
& \le C(T)\left(\frac{\lfloor Nz \rfloor}{N}-\frac{\lfloor Ny \rfloor }{N}\right).
\end{align*}
Since $0<z-y<1$, we deduce that
\begin{align*}
\E\left(\sup_{0\le t\le T}\vert V^{N,y,z}_t\vert^2 \vert \mathcal{G}^{x,y}\right)&\le 3(1+C(T))\left(\frac{\lfloor Nz \rfloor}{N}-\frac{\lfloor Ny \rfloor }{N}\right)\\ 
&+3\beta^2T\int_0^T \E\left(  \sup_{0\le s \le r}(V^{N,y,z}_s)^2\vert  \mathcal{G}^{x,y}\right)dr,
\end{align*}
From this and Gronwall's lemma we deduce that there exists a constant $K_1>0$ such that
\begin{align}
\E\left(\sup_{0\le t\le T}\vert V^{N,y,z}_t\vert^2 \vert \mathcal{G}^{x,y}\right)\le K_1\left(\frac{\lfloor Nz \rfloor}{N}-\frac{\lfloor Ny \rfloor }{N}\right).
\end{align}
Similary we have
$$\E\left[\sup_{0\le t\le T}\left(V^{N,x,y}_s\right)^2\right]\le K_1\left(\frac{\lfloor Ny \rfloor}{N}-\frac{\lfloor Nx \rfloor }{N}\right),$$
Since $0\le y-x<z-x$ and $0\le z-y<z-x$, we deduce that
\begin{align*}
\E\left[\sup_{0\le t\le T}\vert V^{N,x,y}_t\vert^2\times \sup_{0\le t\le T}\vert V^{N,y,z}_t\vert^2\right]\le K_1^2\left(\frac{\lfloor Nz \rfloor}{N}-\frac{\lfloor Nx \rfloor }{N}\right)^2 ,
\end{align*}
hence the result. \\
% \noindent{\sc End of the proof of Proposition \ref{4tightD}}
 
\noindent{\sc Proof of Theorem \ref{4thConv}}
%If the trajectories of $t\to Z^{N,x}_t$ would be continuous for any $x$, then Proposition \ref{4critBI}
%would tell us that the sequence $\{Z^{N,x}_\cdot,x\ge0\}_{N\ge1}$ would be tight in 
%$D([0,\infty);C([0,T]))$, for all $T>0$. \epf
We now show that for any $T>0$,
$$\{Z^{N,x}_t,\ 0\le t\le T,\ x\ge0\}\Rightarrow\{Z^x_t,\ 0\le t\le T,\ x\ge0\}$$
in $D([0,\infty);D([0,T],\R_+))$.
From Theorems 13.1 and 16.8 in \cite{4BI}, since from Proposition \ref{4findimconv}, for all $n\ge1$,
$0<x_1<\cdots< x_n$, 
$$(Z^{N,x_1}_\cdot,\ldots,Z^{N,x_n}_\cdot)\Rightarrow 
(Z^{x_1}_\cdot,\ldots,Z^{x_n}_\cdot)$$ in $D([0,T];\R^n)$,  it suffices to show that
for all $\bar{x}>0$, $\epsilon$, $\eta>0$, there exists $N_0\ge1$ and $\delta>0$ such that for all $N\ge N_0$,
\begin{equation}\label{4w''}
\P(w_{\bar{x},\delta}(Z^N)\ge\epsilon)\le\eta,
\end{equation}
where for a function $(x,t)\to z(x,t)$
$$w_{\bar{x},\delta}(z)=\sup_{0\le x_1\le x\le x_2\le\bar{x},x_2-x_1\le\delta}
\inf\left\{\|z(x,\cdot)-z(x_1,\cdot)\|,\|z(x_2,\cdot)-z(x,\cdot)\|\right\},$$
with the notation $\|z(x,\cdot)\|=\sup_{0\le t\le T}|z(x,t)|$.
But from the proof of Theorem 13.5 in \cite{4BI}, \eqref{4w''} for $Z^N$ follows from Proposition \ref{4tightD}

\section[Generalized Ray Knight Theorem]{Ray Knight representation of a general Feller diffusion}
In this section we establish a Ray-Knight representation of Feller's 
branching diffusion solution of \eqref{4fA}, in terms of the local time of a reflected Brownian motion $H$ with a drift that depends upon the local time accumulated by $H$ at its current level, through the function $f'$ where $f$ is a function satisfying the following hypothesis.\\
\textbf{Hypothesis B:} 
 $f\in C^1(\mathbb{R}_+;\mathbb{R_+})$, $f(0)=0$ and there exist a constant $\beta >0$ such that
\begin{align*} f'(x)\le \beta \quad \forall x\ge0.
\end{align*}
Note that hypothesis \textbf{B} follows from  hypothesis \textbf{A}  if we assume that $f$ is  differentiable. 

The proof we give here  is purely in terms of stochastic analysis, and is inspired by previous work of Norris, Rogers and Williams \cite{4NRW} and Pardoux, Wakolbinger \cite{4PW}. We specify an SDE for a process $\left(H_s\right)$, from which the generalized Feller's diffusion solution of  \eqref{4fA} can be read off from reflected Brownian motion with a drift that is a function of the local time accumulated by $H$ at its current level. 
%The underlying population model is explained by a discrete model exposed in section 1. The process $\left(H_s\right)$ will be reflected Brownian motion with a drift that depends on the local time accumulated at the current level $H_s$ up to time $s$.  
One way to understand the form of the drift is to see $\left(H_s\right)$ as the limit of the exploration process $H^N$ of the forest of random  trees associated to $Z^{N,x}$. 
Precisely,  fix $z\in C(\mathbb{R}_+;\mathbb{R_+})$, the set of continuous functions from $\R_+$ into $\R_+$ and  consider the stochastic differential equation
\begin{equation} \label{4eqHz}
H^z_s=B_s+\frac{1}{2}\int_0^s f'(z(H^z_r)+L^z_r(H^z_r))dr+\frac{1}{2}L^z_s(0),
\end{equation}
where $B$ is a standard Brownian motion, and for $s, t \ge 0$ $L^z_s(t)$ is the  local time accumulated by $H^z$ at level t up to time $s$. For $x>0$ define \[ S_x=\inf \left\lbrace r>0: L^z_s(0)>x\right\rbrace  \text{ and } S=\sup_{x>0}S_x.\]
We first suppose that $f$ satisfies hypothesis \textbf{B} and the following.\\
\textbf{Hypothesis C:} 
\begin{align*}\exists ~a,~ b\in \R :  \forall z\ge 0, \quad \vert f'(z)\vert \le a z +b. \end{align*}
\subsection{Case where $f'$ satisfies hypothesis \textbf{C}.}
In this subsection we suppose that $f$ verifies hypothesis \textbf{C}. We have
 \begin{proposition}\label{4probor}
For any $z\in C(\mathbb{R}_+;\mathbb{R_+})$, equation \eqref{4eqHz} has a unique weak solution.
\end{proposition}
\bpf
Let $H$ denote Brownian motion reflected above 0 defined on a probability space $\left(\Omega,\mathcal{F}, \P \right)$. $H$ solves  the following equation
\[ H_s=B_s+\frac{1}{2}L_s(0), \]
where $B$ is a $\mathcal{F}_s$ standard Brownian motion, and $L$ is the local time of $H$. Let $$M_s:=\frac{1}{2}\int_0^s f'(z(H_r)+L_r(H_r))dB_r \text{   and  } G_s=\exp\left( M_s-\frac{1}{2}\langle M\rangle_s \right).$$
We will show below that $\E(G_s)=1$, for all $s\ge 0$, which is a sufficient condition for $G$ to be a martingale. By application of the Girsanov theorem, there exists a new probability $\tilde{\P^z}$ on $\left(\Omega,\mathcal{F}\right)$ such that
\[ \dfrac{d\tilde{\P^z}}{d\P}\vert_{\mathcal{F}_s}=G_s,  s\ge 0,\]
where $\left( \mathcal{F}_s,s\ge 0\right)$ is the natural filtration of $H$.
Moreover under $\tilde{\P^z}$, $$\tilde{B}^z_s:=B_s-\frac{1}{2}\int_0^sf'(z(H_r)+L_r(H_r))dr, ~~s\ge 0$$ is a standard Brownian motion. The fact that $\E(G_s)=1$ for any $s\ge 0$ follows thanks to assumption $\bf{C}$ from the existence of a constant $c$  such that 
\begin{align} \label{4condg}\sup_{0\le r\le s}\E(\exp(c (L_r(H_r))^2)<\infty .\end{align}
The inequality  \eqref{4condg} is estabilished in \cite{4PW},  see  Lemma 2 and Lemma 3. The uniqueness is also proved in \cite{4PW} and that argument does not make use of hypothesis \textbf{C}. 
\epf

For $K>0$, we now consider Brownian motion reflected in the interval $[0,K]$
\[H^K=B_s+\frac{1}{2}L^K_s(0)-\frac{1}{2}L^K_s(K^-),\]
defined on $\left(\Omega,\mathcal{F},\P \right)$, where $L^K$ denotes the local time of $H^K$. Define for $x>0$
\[ S^K_x=\inf\left\lbrace s>0;L^K_s(0)>x\right\rbrace . \] 
We again define
$$M^K_s:=\frac{1}{2}\int_0^sf'(z(H^K_r)+L^K_r(H^K_r))dB_r \text{   and  } G^K_s=\exp\left( M^K_s-\frac{1}{2}\langle M^K\rangle_t\right).$$
The same argument as above shows that $\E(G^K_s)=1$ for all $s\ge 0$. This implies that there exists a probability $\tilde{\P}^{K,z}$ defined on the  measurable space $\left(\Omega,\mathcal{F}\right) $  under which
\begin{align*} \tilde{B}^z_s=H^K-\frac{1}{2}L^K_s(0)+\frac{1}{2}L^K_s(K^-)-\frac{1}{2}\int_0^sf'(z(H^K_r)+L^K_r(H^K_r))dr,\end{align*}
 is a $\tilde{\P}^{K,z}$-Brownian motion. That is the equation 
 \begin{equation}\label{4eqHK} 
 H^K=\tilde{B}_s+\frac{1}{2}L^K_s(0)-\frac{1}{2}L^K_s(K^-)+\frac{1}{2}\int_0^sf'(z(H^K_r)+L^K_r(H^K_r))dr
 \end{equation}
 admits a weak solution. Uniqueness of the weak solution of \eqref{4eqHK} is obtained in a similar way as concerning  \eqref{4eqHz}. Moreover we have that (see again \cite{4PW}) $$\tilde{\P}^{K,z}\left(S_x^K<\infty\right)=1.$$
We now have the following Ray Knight representation.
\begin{proposition}\label{4RNb}
For any $K>0$ and $z\in C(\mathbb{R}_+;\mathbb{R_+})$,  the law of  $\Big\lbrace L^K_{S^K_x}(t),$\\
$0\le t<K \Big\rbrace$ under $\tilde{\P}^{K,z}$ is the same as the law of  $ \left\lbrace Z^{x,z}_t,0\le t<K\right\rbrace$, 
where $Z^{x,z}$ solves the SDE
\begin{align} \label{4Zz}
dZ^{x,z}_t= \left[ f(Z^{x,z}_t+z(t)) -f(z(t))\right] dt+2\sqrt{Z^{x,z}_t}dW_t,~Z^{x,z}_0=x, 
\end{align}
and $W$ is a standard Brownian motion. 
\end{proposition}
The proof of this Proposition  is similar to that done in \cite{4PW} in the quadratic case. We will give below some details of the proof in the more general case without  hypothesis \textbf{C}.
\subsection{Existence and uniqueness of weak solution of \eqref{4eqHK} without hypothesis \textbf{C}}
Now we do not suppose anymore that $f$ satisfies hypothesis \textbf{C}. We have 
\begin{proposition}\label{4sol}
For any $K>0$, $z\in C(\mathbb{R}_+;\mathbb{R_+})$, there exists a probability $\tilde{\P}^{K,z}$ under which equation \eqref{4eqHK} has a unique weak solution on the random interval $[0,S^K)$,
where $S^K=\sup_{x\ge 0}S^K_x$.
\end{proposition}
\bpf
 Consider again, for $K>0$, the Brownian motion $H^K$ reflected in the interval $[0,K]$ defined on $\left(\Omega, \mathcal{F},\P\right)$. 
\begin{equation*}
H^K=B_s+\frac{1}{2}L^K_s(0)-\frac{1}{2}L^K_s(K^-).
\end{equation*}
For $n\ge 1$,  we define the function $g_n(r)=f'(n\wedge r)$. It is clear that there exist two constants $a,b\ge 0$ such that 
$\vert g_n(r)\vert \le ar+b$. Thanks to Proposition \ref{4probor}, there exits for each $n\ge1$ a probability $\tilde{\P}^{K,z,n}$ such, 
\[ \dfrac{ d\tilde{\P}^{K,z,n} }{d\P}\vert_{\mathcal{F}_s}=\exp\left\lbrace M^{K,n}_t-\frac{1}{2}\langle M^{K,n} \rangle_t\right\rbrace , s\ge 0,\]
where $M^{K,n}_s:=\frac{1}{2}\int_0^s g_n(z(H^K_r)+L^K_r(H^K_r))dB_r$. Under $\tilde{\P}^{K,z,n}$, 
\begin{equation}
\tilde{B}^{z,n}_s=H^K_s-\frac{1}{2}L^K_s(0)-\frac{1}{2}\int_0^s g_n(z(H^K_r)+L^K_r(H^K_r))dr+\frac{1}{2}L^K_s(K^-),\forall s\ge 0,
\end{equation}
is a standard Brownian motian. 
 For $n\ge 1$, we define the stopping time 
\[ T_n=\inf \left\lbrace s>0; \sup_{0\le t<K}\left[z(t)+L^K_s(t)\right]>n \right\rbrace .\]
We need the following result which is a variant of Theorem 1.3.5 from Stroock-Varadhan \cite{4SV}, whose proof is very similar to that in \cite{4SV}.
\begin{theorem} \label{4thr}
Let $\Omega=C(\R_+,\R_+)$ be the canonical path space with its canonical filtration $\left\lbrace \mathcal{F}_t\right\rbrace$, and let $\left( T_n\right) $ be an increasing sequence of stopping times satisfying $T_n\le S^K$ $a.s.$ $\forall n\ge1$. Suppose there is a sequence $\left( \P_n\right)$ of probabilities on $(\Omega,\mathcal{F})$ such that
\begin{itemize}
\item $\P_{n+1}$ agrees with $\P_n$ on $\mathcal{F}_{T_n}$; \\
\item for each  $x>0$, \[\P_n \left( T_n< S^K_x\right)\rightarrow 0  \text{  as  }  n\rightarrow\infty.\] 
\end{itemize}
Then  there exists a probability $\P$ on $(\Omega,\mathcal{F}_{S^K})$ such that for each $n$, \[ \P=\P_{n} \text{   on   } \mathcal{F}_{T_n}. \] 
\end{theorem}
This proves the existence of a probability $\tilde{\P}^{K,z}$ on $(\Omega,\mathcal{F}_{S^K})$, provided we show that for all $x>0$,  $$ \tilde{\P}^{K,z,n}\left(T_n<S^K_x\right)\rightarrow 0 \text{ as } n\rightarrow\infty.$$
We have
\[ \tilde{\P}^{K,z,n}\left( T_n<S^K_x\right)=\tilde{\P}^{K,z,n}\left( \sup_{0\le t<K}L^K_{S^K_x}(t)>n\right).\]
%Futhermore we have that  $g_n$ satisfies hypothesis \textbf{C} with $a=0$ and $b=\sup_{0\le r\le n}\vert f'(r)\vert $. Consequentely, by 
From Propostion \ref{4RNb}, under $\tilde{\P}^{K,z,n}$, $(L^K_{S^K_x}(t),0\le t< K )$ has the same law as  $(Z^{x,z,n}_t,0\le t< K)$ solution of the SDE
\begin{align*}
dZ^{x,z,n}_t=\left(\int_{z(t)}^{z(t)+Z^{x,z,n}_t}g_n(u)du\right)dt+2\sqrt{Z^{x,z,n}_t}dW_t,~~Z^{x,n}_0=x.
\end{align*}
For any $x\ge 0$, consider  the process $\tilde{Z}^x$, which is solution of the SDE 
\begin{align*}
\tilde{Z}^x_t=x+\beta\int_0^t\tilde{Z}^x_rdr+2\int_0^t\sqrt{\tilde{Z}^x_r}dWr.
\end{align*}
By a well known comparison theorem for one dimensional SDEs, see \cite{4RY} theorem $X.3.7$, we have that  for any $x\ge 0$ and $z\in C(\mathbb{R}_+;\mathbb{R_+})$, $Z^{x,z,n}\le \tilde{Z}^x$ $a.s.$ .
We then have that 
\begin{align*}
\tilde{\P}^{K,z,n}\left( T_n<S^K_x\right)&=\tilde{\P}^{K,z,n}\left( \sup_{0\le t<K}L^K_{S^K_x}(t)>n\right) \\
&=\P\left( \sup_{0\le t<K}Z^{x,z,n}_t>n\right) \\
&\le \P\left( \sup_{0\le t<K}\tilde{Z}^x_t>n\right)\rightarrow 0 \text{ as } n\to\infty.
\end{align*}
We thus have proved for all $K>0$ and $z\in C(\mathbb{R}_+;\mathbb{R_+})$, the existence of a probability $\tilde{\P}^{K,z}$ under which on $[0,S^K)$,  
\begin{align*}
\tilde{B}^z_s=H^K_s-\frac{1}{2}L^K_s(0)+\frac{1}{2}\int_0^sf'(z(H^K_r)+L^K_r(H^K_r))dr+\frac{1}{2}L^K_s(K^-)
\end{align*}
is a standard Brownian motion. Uniqueness is obtained in a similar way as in \cite{4PW}. Hence the Proposition.
\epf
For any $z\in C(\R_+;\R_+)$, we have the following Ray Knight representation.
\begin{proposition}\label{p3.5}
For any $K>0$, $z\in C(\R_+;\R_+)$ and $x\ge 0$, the law of $\left( L^K_{S^K_x}(t),0\le t<K\right)$ under $\tilde{\P}^{K,z}$ is the same as the law of  $\left( Z^{x,z}_t,0\le t<K\right)$
\end{proposition}
\bpf 
For $K>0$ and $z\in C(\R_+)$, we work under the probability measure $\tilde{\P}^{K,z}$. Using Tanaka's formula, we have for any $r\ge0$ and $0\le t<K$, the following identity
\begin{align} \label{4tanfor} (H^K_r-t)^-=(-t)^-+\int_0^r{\bf1}_{\lbrace H^K_s\le t\rbrace}dH^K_s+\frac{1}{2}L^K_r(t)  \end{align}
Recall that for any $x\ge 0$, $\P^{K,z}(S^K_x<\infty)=1$. Hence from \eqref{4tanfor},
\[ L^K_{S^K_x}(t)=2\int_0^{S^K_x}{\bf1}_{\lbrace H^K_s\le t\rbrace}dH^K_s .\]
Combining with equation \eqref{4eqHK}, we get
\begin{align*}
L^K_{S^K_x}(t)=x+2\int_0^{S^K_x}{\bf1}_{ \lbrace H^K_s\le t\rbrace}d\tilde{B}_s+\int_0^{S^K_x}{\bf1}_{\lbrace H^K_s\le t\rbrace }f'(z(H^K_s)+L^K_s(H^K_s) )ds.
\end{align*}
From the generalized occupation time formula (see Exercise $VI.1.15$ in \cite{4RY}), we  obtain 
\begin{align*} \int_0^{S^K_x}{\bf1}_{\lbrace H^K_s\le t\rbrace }f'(z(H^K_s)+L^K_s(H^K_s))ds &=\int_0^t\int_0^{S^K_x}f'(z(r)+L^K_s(r))dL^K_s(r)dr\\
&= \int_0^t \left(f(z(r)+L^K_{S^K_x}(r))-f(z(r))\right)dr .
\end{align*}
The key idea of the  proof is now to introduce the ``excursion filtration", as  in \cite{4NRW} and \cite{4PW}.
For any $0\le t<K$ and $s\ge0$, let define
\begin{align*}
A_s(t)&:=\int_0^s {\bf1}_{\lbrace H^K_r\le t\rbrace}dr, ~~\tau(r,t)=\inf\left\lbrace s>0; A_s(t)>r\right\rbrace ,\\
J(s,t)&:=\int_0^s {\bf1}_{\lbrace H^K_r\le t\rbrace}d\tilde{B}_r~~~~\xi(r,t):=J(\tau(r,t),t)\\
\mathcal{F}_{(r,t)}&:=\sigma(\xi(r,t):0\le u\le r),~~ \varsigma_t=\mathcal{F}_{(\infty,t)},\\
N_t &:=\int_0^{S^K_x}{\bf1}_{\lbrace H^K_s\le t\rbrace}d\tilde{B}_s.
\end{align*}
For fixed $t$, the process $\left( J(s,t),s\ge 0\right) $ is a martingale with respect to $\mathcal{F}_s$, while $\left( \xi(r,t),r\ge 0\right)$ is a $\mathcal{F}_{(r,t)}$-martingale and its quadratric variation equals $r$. Consequently 
$\left( \xi(r,t),r\ge 0\right)$ is a $\mathcal{F}(r,t)$-Brownian motion. We then have the
\begin{lemma}
The process $\left( N_t,0\le t<K\right)$ is a continuous $\varsigma_t$-martingale with its quadratic variation given by
\[ \langle N \rangle_t=4\int_0^t L^K_{S^K_x}(r)dr.\]
\end{lemma}
This result is an easy consequence of Theorem 1 in \cite{4NRW}. By the martingale representation theorem, we deduce that there exists a Brownian motion $W$ such that 
$$N_t=2\int_0^t\sqrt{L^K_{S^K_x}(r)}dW_r.$$
Consequently for all $0\le t< K$,  
\[ L^K_{S^K_x}(t)= x+\int_0^t\left( f(z(r)+L^K_{S^K_x}(r))-f(z(r))\right)dr+2\int_0^t\sqrt{L^K_{S^K_x}(r)}dW_r.\]
\epf 
\subsection{Ray Knight theorem in the subcritical case} 
We first prove the following proposition (recall the definition \eqref{lambda} of $\Lambda(f)$)
\begin{proposition}
Suppose that  $f$ satisfies hypothesis \textbf{B} and 
$\Lambda(f) =\infty$.  Then equation \eqref{4eqHz} admits a unique weak solution on $[0,S)$.
\end{proposition}
\bpf
For $x>0$ and $K>0$, let define
\[ \Omega ^{K,x}= \left\lbrace  \sup_{[0,S^K_x]} H^{K'}<K,~~\forall~ K'>K \right\rbrace .\]
For any $x\ge 0$ and $z\in C$, since we are in the subcritical case, there exists $T_{x,z}<\infty$ $a.s.$ such that , $Z^{x,z}_t=0$, $\forall$ $t\ge T_{x,z}$. We deduce from this and Proposition \ref{p3.5} that for any  fixed $x\ge 0$,
\[ \Omega=\cup_{K>0}\Omega^{K,x} \text{ a.s. } .\]
Note that the family of events $\left\lbrace \Omega^{K,x},~K>0\right\rbrace $ is increasing, and on $\Omega^{K,x}$,  $H^{K'}=H$ a.s., for any $K' >K$. We can define a probability $\tilde{\P}^{z,x}$ on $\left(\Omega, \mathcal{F}_{S_x}\right)$ such that $\tilde{\P}^{z,x}$=$\tilde{\P}^{K,z}$ on $\Omega^{K,x}$. Under $\tilde{\P}^{z,x}$, on $[0,S_x]$
\[\tilde{B}^z_s=H_s-\frac{1}{2}\int_0^s f'(z(r)+L_r(H_r))dr-\frac{1}{2}L_s(0) \]
is a standard Brownian motion. This proves that \eqref{4eqHz} has a  weak solution on $[0,S_x]$ whose uniqueness  can be proved as in \cite{4PW}.
ON A A NOUVEAU BESOIN DU TH DE SV !
We can deduce that there exists a probabilty $\tilde{\P}^z$ under which, on  $[0,S)$ $\tilde{B}^z$ is a standard Brownian motion, where \[ S=\sup_{x\ge 0}S_x . \] 
\epf 
For $z\equiv 0$, we write $\tilde{\P}=\tilde{\P}^{z}$. The following statement is a generalized Ray Knight theorem in the subcritical case.
\begin{theorem} \label{4RN} Suppose that  $f$ satisfies Hypothesis \textbf{B} and 
$\Lambda(f) =\infty$. Then the law of the random fields  $\left\lbrace L_{S_x}(t), t\ge 0, x\ge 0\right\rbrace$ under the probability $\tilde{\P}$  is the same as the law of  $\left\lbrace Z^x_t,t\ge 0, x\ge 0\right\rbrace $. 
\end{theorem}
We first establish the following Proposition.
\begin{proposition}\label{4RNx}
Assume that the two assumptions of Theorem \ref{4RN} holds. Then for any $x$ and $z\in C(\R_+;\R_+)$ fixed, 
the law of $\left \lbrace L_{S_x}(t), t\ge 0\right\rbrace$ under $\tilde{\P}^z$ coincides with is the law of  $\left\lbrace Z^{x,z}_t,t\ge 0\right\rbrace$.
\end{proposition}
\bpf
We have that for any $K>0$, $z\in C(\R_+;\R_+)$, under $\tilde{\P}^{K,z}$, $\Big( L^K_{S^K_x}(t),$ \\
$0\le t<K\Big)$ has the same law as $\left( Z^{x,z}_t,0\le t<K\right)$. A consequence of this  is that for any $0<K<K',$
\begin{align}\label{4idenlaw}
\left\lbrace L^K_{S^K_x}(t),0\le t<K \right\rbrace \stackrel{(d)}=\left\lbrace L^{K'}_{S^{K'}_x}(t),0\le t<K  \right\rbrace. 
\end{align}
It now follows that for any $K$, under $\tilde{\P}^z$, $\left( L_{S_x}(t), 0\le t<K \right)$ has the same law as $\left( L^K_{S^K_x}(t), 0\le t<K \right)$ under $\tilde{\P}^{K,z,x}$. We then obtain that for any $K>0$ 
\[ \left( L_{S_x}(t),0\le t<K \right) \stackrel{(d)}=\left( Z^{x,z}_t,0\le t<K \right).\]
Hence the proposition, letting $K$ go to $\infty$.
\epf
In particular, for $x$ fixed, the law of $\left \lbrace L_{S_x}(t), t\ge 0\right\rbrace$ under $\tilde{\P}$ is the same as the law of $\left\lbrace Z^x_t,t\ge 0\right\rbrace$.
\begin{remark}
The identity  \eqref{4idenlaw} could also obtained from a generalization of  Lemma 2.1 in Delmas \cite{4DJ}. For $0<a<b$, we define the application $\pi^{a,b}$ with maps continuous trajectories with value in $[0,b]$ into trajectories with values in $[0,a]$ as follows. If $u\in C(\R_+,[0,b])$,
\[ \rho_u(s)=\int_0^s{\bf1}_{u(r)<a}dr,~~~~~~\pi^{a,b}(u)(s)=u(\rho^{-1}_u(s)).\]
The following equality in law holds.
\[\pi^{a,b}(H^b)\stackrel{(d)}=H^a.\]
This identity together with the strong Markov property of the Brownian motion implies \eqref{4idenlaw}.
\end{remark}
%\subsection*{Completion of the proof of the Theorem \ref{4RN}}
\noindent{\sc Proof of Theorem \ref{4RN}}
Recall that $\left( Z^x_., x\ge 0\right) $ is a Markov process with value in the space of continuous paths from $\R_+$ into $\R_+$ with compact support. From Proposition \ref{4RNx} with $z\equiv 0$, its marginal laws coincide with those of $L_{S_x}(.)$. We now check that   $\left(L_{S_x}(.),x\ge 0\right)$ is a Markov process. This follows readily from the fact that for any $0\le x<y $, conditionnaly upon $\left(L_{S_{x'}}(.),x' \le x\right)$ and given $L_{S_x}(.)=z(.)$, on $[0,S_y]$ the process $H^x_s:=H_{S_x+s}$ solves the SDE
\begin{align*}
H^x_s=\bar{B}_s+\frac{1}{2}\int_0^s\left(f'(z(H^x_r)+L^z_r(H^x_r))\right)dr+\frac{1}{2}L^z_s(0),
\end{align*}
where $\bar{B}$ is a Brownian motion independent of $\left( L_{S_{x'}}(t),x'\le x,0\le t\le S_x\right)$ and $L^z$ denotes the local time of $H^x$, which is also the additional local time accumulated by $H$ after time $S_x$. To complete the proof of the theorem it now suffices to prove that  for any $x,y\ge 0$ the conditional law of $\left( L_{S_{x+y}}(t),t\ge 0\right)$ given $\left( L_{S_x}(t),t\ge0\right)$ is the same as the conditional law of $\left( Z^{x+y}_t,t\ge\right)$ given $\left( Z^x_t,t\ge 0\right)$. Conditioned upon $L_{S_x}(.)=z(.)$, $L_{S_{x+y}}(.)-L_{S_x}()$ is the collection of local times accumulated by $H^x$ up to time $S_{y}$, and it has the same law as $L^z_{S_y}(.)$ while conditionally upon $Z^x=z(.)$, $Z^{x+y}-Z^x$ has the same law as  $Z^{y,z}$. The identity of those two laws has been established in Proposition \ref{4RNx}.
% Thus the two objects have the same transitions probabilities. As consequence of that, the assertion of the Theorem follows from Proposition \ref{4RNx}.
\epf 
We can deduce from the Proposition \ref{4RNx} and the occupation time formula that 
\begin{corollary} Suppose that  $f$ satisfies Hypothesis \textbf{B} and $\Lambda(f) =\infty$. We have  $$\forall x\ge0, \quad \tilde{\P}\left( S_x<\infty\right)=1. $$
\end{corollary}
\bpf
Let $g(h)=1$, for any $h>0$. By the occupation times's formula, we have
\begin{align*}
S_x& =\int_0^{S_x}g(H_r)dr\\
&=\int_0^{\infty}L_{S_x}(t)dt=\int_0^{T_0^x}Z^x_rdr<\infty \text{ a.s.} 
\end{align*}
\epf
Note that $S_x$ is the total mass  of the process $\left(Z^x_t,t\ge0\right)$.
{}    

\end{document}